\input amstex
\documentstyle{amsppt}
%\magnification=1200
\pagewidth{5.5in}
%% STANDARD ABBREVIATIONS
%\magnification=1200

\let\bs\backslash

\let\p\prime    
\let\ol\overline
\let\ul\underline
\define\a{{\alpha}}

\define\K{{\Cal K}}

\define\fL{{\Cal L}}

\define\X{{\Cal X}}
\define\Y{{\Cal Y}}
\define\CC{{\Bbb C}}
\define\QQ{{\Bbb Q}}
\define\Q{{\Bbb Q}}

\define\PP{{\Bbb P}}

\define\MHS{\text{\rm MHS}}
\define\Spec{\text{\rm Spec}}
\define\Ext{\text{\rm Ext}}
\define\ext{\text{\rm Ext}^1_{\text{\rm MHS}}}
\define\CH{\text{\rm CH}}
\define\ch{\text{\rm CH}}

\define\alg{\text{\rm alg}}

\define\IM{\text{\rm Im}}
\define\im{\text{\rm im}}
\define\codim{\text{\rm codim}}
\define\coker{\text{\rm coker}}
\define\Griff{\text{\rm Griff}}

\define\cl{\text{\rm cl}}
\define\level{\text{\rm Level}}

\define\tY{{\tilde Y}}
\define\ds{\displaystyle}
\define\dt{\bullet}

\topmatter
\title   Beilinson's Hodge Conjecture for K$_{1}$ Revisited \endtitle
\leftheadtext{S-J. Kang \& J. D. Lewis}
\rightheadtext{Beilinson-Hodge Conjecture}
\author Su-Jeong Kang {\rm and} James D. Lewis  \endauthor
\address  Department of Mathematical and Statistical Sciences, 
University of Alberta, Edmonton, Alberta T6G 2G1, Canada
\endaddress
\email sjkang\@math.ualberta.ca, lewisjd\@ualberta.ca \endemail
\thanks First author supported by a PIMS postdoctoral fellowship.
Second author is partially supported by a grant 
from the Natural 
Sciences and Engineering Research Council of Canada.\endthanks
\abstract
Let $U/\CC$ be a smooth  quasiprojective variety
and $\CH^{r}(U,1)$
a special instance of Bloch's higher Chow groups
(\cite{Blo}). Jannsen was the
first to show that the cycle class map $\cl_{r,1} : 
\CH^{r}(U,1)\otimes \QQ \to 
\hom_{\MHS}\big(\QQ(0),H^{2r-1}(U,\QQ(r))\big)$ is not in general 
surjective, contradicting an earlier conjecture of Beilinson. In this
paper we give a refinement of Jannsen's counterexample, and 
further show that the aforementioned cycle class map becomes
surjective at the generic point.
\endabstract
\endtopmatter

\document
\head \S0. Introduction \endhead

Let $U/\CC$ be a smooth quasiprojective variety, 
$\CH^{r}(U,m;\QQ) := \CH^{r}(U,m)\otimes\QQ$, where
$\CH^{r}(U,m)$ is Bloch's higher Chow group, and 
$$
\cl^{U}_{r,m} : \CH^{r}(U,m;\QQ) \to \Gamma\big(H^{2r-m}(U,\QQ(r))\big)
:= \hom_{\MHS}\big(\QQ(0),H^{2r-m}(U,\QQ(r))\big),
$$
the Betti cycle class map. If $m=0$, then the Hodge conjecture
(classical form) implies that $\cl^{U}_{r,0}$ is surjective. Beilinson
(\cite{Be}) once conjectured that $\cl^{U}_{r,m}$ is 
always surjective. It was Jannsen (\cite{J3}) who was the
first to find a counterexample, in the case $m=1$, where the complex
numbers $\CC$ are used in an essential way. In contrast to this,
one expects the surjectivity of $\cl^{U}_{r,m}$ in the case where $U$
is obtained via base extension from a variety defined over
a number field. It turns out that Jannsen's counterexample
is indeed very special, being the complement of a closed subscheme
of codimension $r$ in a projective algebraic manifold $X$. 
In contrast to this, we show rather easily  that
the corresponding limit map
$$
\cl_{r,1} : \CH^{r}(\Spec(\CC(X)),1;\QQ) \to \Gamma\big(H^{2r-1}(
\CC(X),\QQ(r))\big),
$$
is onto, provided that a certain reasonable conjectural
type statement in \cite{J1} holds, and unconditionally in the case
$r = \dim X$. We review  this in Example 2.7 below. 
Broadly speaking, besides providing a closer examination
of this cycle class map, we want to also consider the
relative situation as well.
Consider a  proper morphism $\rho :  \X \to S$ of
smooth complex quasiprojective varieties. Let $\eta$ be 
the generic point of $S$.
\proclaim{Question 0.1} 
{\rm Is the induced map
$$
\cl_{r,m}^{\eta} : \CH^r(\X_{\eta},m;\QQ) \to \Gamma\big(
H^{2r-m}(\X_{\eta},\QQ(r))\big),
$$
surjective?  
[Here $H^{2r-m}(\X_{\eta},\QQ(r)) = 
\lim_{{\buildrel\to\over {U\subset S}}}H^{2r-m}(\rho^{-1}(U)(\CC),
\QQ(r)).]$}
\endproclaim
Note that 0.1 in the case $m=0$ is equivalent to
the classical Hodge conjecture, using a standard localization
argument (\cite{J3}).
Based on our limited results in \S2 for the case $m=1$, 
we feel that the answer to this
question is yes. The main reason for considering
the relative situation is as follows. 
As a formal consequence of M. Saito's theory
of mixed Hodge modules (see \cite{As}, \cite{KL}), there is an exact sequence:
$$
0\to
{\Ext}^1_{\text{\rm PMHS}}\big(\QQ(0),H^{\nu-1}
(\eta,R^{2r-\nu-m}\rho_{\ast}\QQ(r))\big)\to
\left\{\matrix \text{\rm Germs\ of}\\ 
\text{\rm higher\ order}\\
\text{\rm normal\ functions}\\
\text{\rm at}\ \eta
\endmatrix\right\}
$$
$$
\to\hom_{\MHS}\big(\QQ(0),H^{\nu}({\eta},R^{2r-m-\nu}\rho_{\ast}\QQ(r))
\big)\to 0
$$
Here PMHS stands for the category of graded polarizable MHS 
(mixed Hodge structures). The key point is, is there 
lurking a generalized Poincar\'e
existence theorem for higher normal functions? (Namely,
are these normal functions cycle induced?) An affirmative answer 
to 0.1, would imply a generalized Poincar\'e
existence theorem, and conversely such an existence theorem should
lead to an affirmative answer to 0.1.
\bigskip
In \S2 we discuss a relative version of the cycle
class map on ``the generic fiber'', as
well as examples pointing to an affirmative answer to 0.1
in the case $m=1$. In \S3 we study the relative
situation more carefully, by looking at those open subsets
of a complete base variety $\ol{S}$, obtained by deleting codimension
$p$ subvarieties from $\ol{S}$. The cokernel of a resulting
limiting cycle class map admits a rather explicit description 
in terms of two filtrations on a given Chow group. The main results
of \S3 are summarized in Corollary 3.6.
For the remainder of the paper, we restrict to the case
$\X = S$, with $\rho =$ identity, and where $X := \ol{\X}$ 
is a projective algebraic manifold. We then describe the aforementioned
cokernel in terms of the complexity of the Chow ring
$\CH^{\ast}(X;\QQ)$. This is accomplished
with the aid of our main technical result in Theorem 4.2,
with subsequence corollary consequences given in \S5. Finally
in \S6 we present some key examples.

\head \S1. Notation\endhead

The following notation will be used throughout the paper.
Let $X$ be a smooth complex projective variety of dimension $d$.
\medbreak\noindent
(i) $\Gamma ( - ) = \hom_{\MHS}(\Q(0), -)$
\medbreak\noindent
(ii) $J( - ) = \ext(\Q(0), - )$
\medbreak\noindent
(iii) $J^r(X) = \ext\big(\Q(0), H^{2r-1}(X,\Q(r))\big)$
\medbreak\noindent
(iv) $\ch^r(X;\Q) = \ch_{d-r}(X;\Q)$ is the Chow group of cycles 
of codimension $r$ ($\dim d-r$) in $X$ modulo rational equivalence,
tensored with $\Q$,
$\ch_{\alg}^r(X;\Q) \subset \ch_{\hom}^r(X;\Q)\subset \ch^r(X;\Q)$
are the subgroups of cycles algebraically equivalent to zero (resp.
nullhomologous), and $\text{Griff}^r(X)_{\Q} = 
\ch^r_{\hom}(X;\Q) / \ch^r_{\alg}(X;\Q)$ is 
the Griffiths group ($\otimes \Q$). More generally,
one has the higher Chow groups $\CH^{r}(X,m)$ due to Bloch, and defined
in \cite{Blo}. 
\medbreak\noindent
(v) $AJ : \ch^r_{\hom}(X;\Q) \to J^r(X)$ the Abel-Jacobi map
\medbreak\noindent
(vi) $\ch^r(X;\Q)_{AJ} = \ker(AJ)$
\medbreak\noindent
(vii) HC = Hodge conjecture (classical form)
\medbreak\noindent
(viii) GHC = (Grothendieck amended) general Hodge conjecture 
(\cite{G}, \cite{L}).

\head \S2. The cycle map on the generic fiber\endhead

The setting is the following diagram
$$
\matrix
\X&\hookrightarrow&\ol{\X}\\
&\\
\rho\biggl\downarrow\ &&\ \ \biggr\downarrow\ol{\rho}\\
&\\
S&\hookrightarrow&\ol{S}
\endmatrix\tag{2.1}
$$
where $\ol{\X}$ and $\ol{S}$ are nonsingular complex 
projective varieties, $\ol{\rho}$ is a dominating flat morphism,
$D\subset \ol{S}$ a divisor, $\Y:= \ol{\rho}^{-1}(D)$,
$S := \ol{S}\bs D$, $\X := \ol{\X}\bs \Y$ and $\rho := \ol{\rho}\big|_{\X}$.
For notational simplicity, we will identify 
$H_{\Y}^{2r-1}(\ol{\X},\QQ(r))$ with its image in $H^{2r-1}(\ol{\X},\QQ(r))$.
There is a short exact sequence
$$
 0\to \frac{H^{2r-1}(\ol{\X},\QQ(r))}{H_{\Y}^{2r-1}(\ol{\X},\QQ(r))} 
 \to H^{2r-1}(\X,\QQ(r)) 
 \to  H^{2r}_{\Y}(\ol{\X},\QQ(r))^{\circ}\to 0,\tag{2.2}
$$ 
and corresponding diagram:
$$
\matrix
\ch^{r}(\X,1;\Q)&\to&\ch^{r}_{\Y}(\ol{\X};\Q)^{\circ}&@>\alpha_{\Y}>>&
\ch^{r}_{\hom}(\ol{\X};\Q)\\
&\\
\cl^{\X}_{r,1}\biggl\downarrow\quad&&\quad \biggl\downarrow\beta_{\Y}&&\quad
\biggl\downarrow \ul{AJ}^{\ol{\X}}\\
&\\
\Gamma\big( H^{2r-1}(\X,\Q(r))\big)&\hookrightarrow&\Gamma\big(
H^{2r}_{\Y}(\ol{\X},\Q(r))^{\circ}\big)&\to&J\biggl(
\frac{H^{2r-1}(\ol{\X},\QQ(r))}{H_{\Y}^{2r-1}(\ol{\X},\QQ(r))}\biggr)
\endmatrix\tag{2.3}
$$
where $\ul{AJ}^{\ol{\X}}$ is a corresponding reduced Abel-Jacobi map. Let us
assume that $\beta_{\Y}$ is surjective (such is the case 
if the Hodge conjecture holds.) If we 
apply the snake lemma, we arrive at 
$$
\coker(\cl^{\X}_{r,1}) \simeq 
\frac{\ker\biggl[ \ul{AJ}^{\ol{\X}}\big|_{\IM(\alpha_{\Y})}:
\IM(\alpha_{\Y}) \to J\biggl(\frac{H^{2r-1}(
\ol{\X},\Q(r))}{H_{\Y}^{2r-1}(\ol{\X},\Q(r))}\biggr)
\biggr]}{\alpha_{\Y}\big(\ker(\beta_{\Y})\big)}. \tag{2.4}
$$
Now take the limit over all $D\subset \ol{S}$ to arrive at an
induced cycle map:
$$
\cl^{\eta}_{r,1} : \CH^{r}(\X_{\eta},1;\QQ) \to
\Gamma\big(H^{2r-1}(\X_{\eta},\QQ(r))\big). \tag{2.4.1}
$$
where $\eta$ is the generic point of $\ol{S}$. We arrive at:
$$
\frac{\Gamma\big(H^{2r-1}(\X_{\eta},\QQ(r))\big)}{\cl^{\eta}_{r,1}\big(
\CH^{r}(\X_{\eta},1;\QQ)\big)} \simeq\
\frac{\ker\biggl[\K
@>\ul{AJ}>> J\biggl(
\frac{H^{2r-1}(\ol{\X},\QQ(r))}{N_{\ol{S}}^{1}H^{2r-1}(
\ol{\X},\QQ(r))}\biggr)\biggr]}{N_{\ol{S}}^{1}\CH^{r}(\ol{\X};\QQ)},\tag{2.5}
$$
where $\K := \ker [\CH^{r}_{\hom}(\ol{\X};\QQ)\to
\CH^{r}(\X_{\eta};\QQ)]$, and
$N_{\ol{S}}^{q}\CH^{r}(\ol{\X}) \subset 
\CH_{\hom}^{r}(\ol{\X};\Q)$ is the subgroup
generated by cycles which are homologous to zero on some
codimension $q$ subscheme of $\ol{\X}$ obtained from a  (pure) codimension
$q$ subscheme of $\ol{S}$ via $\ol{\rho}^{-1}$, and where  
$N_{\ol{S}}^{q}H^{2r-1}(\ol{\X},\QQ(r))$ is the subspace of the
coniveau $N^qH^{2r-1}(\ol{\X},\QQ(r))$ arising from $q$
codimensional subschemes
of $\ol{S}$ via $\ol{\rho}^{-1}$. In \S 3 we show that
under the assumption of the HC,  (2.5) becomes:
$$
\frac{\Gamma\big(H^{2r-1}(\X_{\eta},\QQ(r))\big)}{\cl^{\eta}_{r,1}\big(
\CH^{r}(\X_{\eta},1;\QQ)\big)} \simeq\
\frac{N_{\ol{S}}^{1}\CH^{r}(\ol{\X};\QQ)\ +\ \ker\big[\K@>{AJ}>> J\big(
H^{2r-1}(\ol{\X},\QQ(r))\big)\big]}{N_{\ol{S}}^{1}\CH^{r}(\ol{\X};\QQ)}\tag{2.6}
$$

\noindent
\underbar{Example 2.7}. \ Suppose that $\ol{\X} = \ol{S}$ with $\ol{\rho}$ 
the identity. In this case (2.6) becomes:
$$
\frac{\Gamma\big(H^{2r-1}(\CC(\ol{\X}),\QQ(r))\big)}{\cl_{r,1}\big(
\CH^{r}(\Spec(\CC(\ol{\X})),1;\QQ)\big)}\  \simeq\
\frac{N^{1}\CH^{r}(\ol{\X};\QQ)\ +\  
\CH^r(\ol{\X};\QQ)_{AJ}}{N^{1}\CH^{r}(\ol{\X};\QQ)},
$$
where $N^1\CH^r(\ol{\X};\QQ)$ is the subgroup of cycles, that are homologous to
zero on codimension $1$ subschemes of $\ol{\X}$. According to Jannsen
(\cite{J1}, p. 227), there is reason to believe that the right
hand side of  2.7 should be conjecturally zero. In particular, 
since $\Spec(\CC(\ol{\X}))$ is
a point, this implies that $\Gamma\big(H^{2r-1}(\CC(\ol{\X}),\QQ(r))\big) 
= 0$ for $r>1$. The reader can easily check that 
$$
\cl_{r,1}\big(\CH^{r}(\Spec(\CC(\ol{\X})),1;\QQ)\big) = 
\Gamma\big(H^{2r-1}(\CC(\ol{\X}),\QQ(r))\big),
$$
holds unconditionally in the case $r=\dim \ol{\X}$.

\bigskip
\noindent
\underbar{Example 2.8}. Here we give some evidence that
the RHS (hence LHS) of (2.6) is zero. 
Suppose $\ol{\X} = X\times \ol{S}$, and let us
assume the condition 
$$
\CH^{r}(\ol{\X};\QQ) = \bigoplus_{\ell=0}^r\CH^{r-\ell}(X;\QQ)\otimes \CH^{\ell}
(\ol{S};\QQ).
$$
An example
situation is when $\ol{S}$ is a flag variety, 
such as a projective space; however conjecturally
speaking, this condition is expected to hold for
a much broader class of examples (see \cite{CL}).
Thus 
$$
\xi \in \CH^{r}(\ol{\X};\QQ) \quad \Rightarrow \quad \xi = \sum_{\ell=0}^r \xi_\ell \in 
\bigoplus_{\ell=0}^r\CH^{r-\ell}(X;\QQ)\otimes \CH^{\ell}(\ol{S};\QQ).
$$
Fix an $\ell$ and write
$$
\xi_{\ell} = \sum_{i=1,j=1}^{N_{\ell},M_{\ell}}\gamma^{\ell}_i\otimes \beta^{\ell}_j.
$$
We can assume that 
$\{[\gamma^{\ell}_1],\ldots,[\gamma^{\ell}_{N_{1,\ell}}]\}$ is a basis for $\QQ[\gamma^{\ell}_1]+\cdots
+\QQ[\gamma^{\ell}_{N_{\ell}}]\subset H^{2r-2\ell}(X,\QQ)$, and that $\{[\beta^{\ell}_1],\ldots,[\beta^{\ell}_{M_{1,\ell}}]\}$ 
is a basis for $\QQ[\beta^{\ell}_1]+\cdots
+\QQ[\beta^{\ell}_{M_{\ell}}]\subset H^{2\ell}(\ol{S},\QQ)$. Therefore we can write
$$
\xi_{\ell} = \sum_{i=1,j=1}^{N_{1,\ell},M_{1,\ell}}\gamma^{\ell}_i\otimes \beta^{\ell}_j \ +\
\sum_{i=1,j=1}^{N_{1,\ell},M_{1,\ell}}\gamma^{\ell}_i\otimes (\beta^{\ell}_j)^{\p} \ +\ 
\sum_{i=1,j=1}^{N_{1,\ell},M_{1,\ell}}(\gamma^{\ell}_i)^{\p} \otimes \beta^{\ell}_j \ +\
\xi^{\p}_{\ell},
$$
where $\xi^{\p}_{\ell} \in \CH^{r-\ell}_{\hom}(X;\QQ)\otimes \CH_{\hom}^{\ell}(\ol{S};\QQ)$,
$(\gamma^{\ell}_i)^{\p}\in \CH^{r-\ell}_{\hom}(X;\QQ)$, $(\beta^{\ell}_j)^{\p}\in \CH_{\hom}^{\ell}(\ol{S};\QQ)$.
It is obvious that $\xi\in \CH_{\hom}^r(\ol{\X};\QQ) \Rightarrow 
\sum_{i=1,j=1}^{N_{1,\ell},M_{1,\ell}}\gamma^{\ell}_i\otimes \beta^{\ell}_j  = 0$ for each $\ell$.  Furthermore,
$AJ(\xi^{\p}) = 0$, where $\xi^{\p} = \sum^{r}_{\ell=0} \xi^{\p}_{\ell}$. This is easily seen from the description
$$
\Ext^1_{\MHS}\big(\QQ(0),H^{2r-1}(\ol{\X},\QQ(r)\big) = \bigoplus_{\ell=0}^r
\left[\matrix \Ext^1_{\MHS}\big(\QQ(0),H^{2r-2\ell-1}(X,\QQ) \otimes 
H^{2\ell}(\ol{S},\QQ))(r)\big)\\
\bigoplus\\
\Ext^1_{\MHS}\big(\QQ(0),H^{2r-2\ell}(X,\QQ) \otimes 
H^{2\ell-1}(\ol{S},\QQ))(r)\big)
\endmatrix\right].
$$
Note that for $\ell \geq 1$, 
$ \sum_{i=1,j=1}^{N_{1\ell},M_{1\ell}}(\gamma^{\ell}_i)^{\p} \otimes \beta^{\ell}_j \ +\
\xi^{\p}_{\ell} \in N^1_{\ol{S}}\CH^r(\ol{\X};\QQ)$.
Recall that $\K =\ker [\CH^{r}_{\hom}(\ol{\X};\QQ)\to
\CH^{r}(\X_{\eta};\QQ)$]. Then:
$$
\K \subset \bigoplus_{\ell\geq 1}\CH^{r-\ell}(X;\QQ)\otimes \CH^{\ell}
(\ol{S};\QQ),
$$
and
$$
\xi\in \K \cap \CH^{r}(\ol{\X};\QQ)_{AJ} \Rightarrow 
AJ_{\ol{S}}\big((\beta^{\ell}_j)^{\p}\big) = 0.
$$
Thus if we assume that  $\ker(AJ_{\ol{S}}) \subset N^1\CH^{\ell}(\ol{S};\QQ)$, 
then we arrive
at $\xi \in N^1_{\ol{S}}\CH^r(\ol{\X};\QQ)$.

\proclaim{Conjecture 2.9} The map $\cl_{r,1}^{\eta} $in (2.4.1) is surjective.
\endproclaim

%############## New section 
%##############################################

\head \S3. A (Relative) Filtration on Chow Groups \endhead

We recall the setting in (2.1) and the definitions of $\K$, 
$N^q_{\ol{S}}\CH^r(\ol{\X};\QQ)$,
$N^q_{\ol{S}}H^{2r-1}(\ol{\X},\QQ(r))$ in (2.5) above,
and introduce the following descending filtration
$\{\fL^p_{\ol{\X}/\ol{S}}\CH^r(\ol{\X};\QQ)\}_{p\geq 0}$:
$$
\fL^{p}_{\ol{\X}/\ol{S}}\ch^r (\ol{\X};\Q) = 
\left.\left\{
\aligned
&\ker\left[\ul{AJ}^p \big|_{\K_p}: \K_p \to J \left( 
\frac{H^{2r-1}(\ol{\X},\Q(r))}{N_{\ol{S}}^{p}H^{2r-1}(\ol{S},\Q(r))} 
\right) \right]\\
&0 
\endaligned\right.\quad
\aligned
&\\
& \text{for $0 \leq p \leq r$}\\
&\\
& \text{for $p \geq r+1$}
\endaligned\right.\tag{3.1}
$$
where $\K_p = \ds{\ker \left[\ch^r_{\hom}(\ol{\X};\Q) \to  
\lim_{{\buildrel\to\over{D \in Z^p(\ol{S})}}} \ch^r(\ol{\X} 
\setminus \big|\ol{\rho}^{-1}(D)\big|;\Q) \right]}$, $Z^p(\ol{S})$ 
is the free abelian group generated by codimension $p$
subvarieties of $\ol{S}$.

\bigskip
\noindent
\underbar{Example 3.2}. Suppose that $X := \ol{\X} = \ol{S}$ where
$\ol{\rho}$ is the identity. In this case put
$\{ \fL^{p}\ch^r(X;\Q)\}_{p \geq 0} = 
\{\fL^p_{\ol{\X}/\ol{\X}}\CH^r(\ol{\X};\QQ)\}_{p\geq 0}$.
We then have:
$$
\fL^{p}\ch^r (X;\Q) = 
\left.\left\{
\aligned
&\ker\left[\underline{AJ}^{p}: \ch^r_{\text{hom}}(X;\Q) \to J \left( 
\frac{H^{2r-1}(X,\Q(r))}{N^{p}H^{2r-1}(X,\Q(r))} \right) \right]\\
&0 
\endaligned\right.\quad
\aligned
& \text{for $0 \leq p \leq r$}\\
& \text{for $p \geq r+1$}
\endaligned\right.
$$
where $\ul{AJ}^{p}: \ch^r_{\text{hom}}(X;\Q) \to J 
\left(\frac{H^{2r-1}(X,\Q(r))}{N^{p}H^{2r-1}(X,\Q(r))} \right)$ is 
the reduced Abel-Jacobi map. Note that $N^r H^{2r-1}(X;\Q)=0 $ and 
hence $\ul{AJ}^r=AJ$ with:
$$
\fL^r \ch^r(X;\Q) = \ker[AJ:\ch^r_{\text{hom}}(X;\Q) \to 
J^r(X)]=\ch^r(X;\Q)_{AJ}
$$
The coniveau filtration $\{N^{p}\ch^r(X;\Q)\}_{p \geq 0}$  on 
Chow groups introduced above appears in \cite{J1}. Recall its 
$p-$th level $N^{p} \ch^r(X;\Q)$ is defined 
to be the subgroup generated by cycles which are homologous to zero 
on some $p-$codimensional possibly reducible subvariety of $X$. 
\bigskip
We  have the following relations:

\proclaim{Proposition 3.3}
For each $p$, we have
$$
N_{\ol{S}}^p \ch^r(\ol{\X};\Q) \subseteq 
\fL^p_{\ol{\X}/\ol{S}} \ch^r(\ol{\X};\Q) 
\subseteq \ch^r_{\hom}(\ol{\X};\Q)
$$
\endproclaim

\demo{Proof}  The first inclusion is a consequence of 
functoriality of the Abel-Jacobi map,
and the second inclusion is by definition. \qed
\enddemo

For $r\geq 1$ we consider 
$$
\cl_{r,1}(\Y) : \ch^r(\ol{\X}\bs \Y,1;\Q )\to 
\Gamma\big(H^{2r-1}(\ol{\X}\bs \Y,\Q(r))\big)
$$
for a closed subscheme $\Y$ of $\ol{\X}$. Set
$$
\cl^p_{r,1} := \lim_{{\buildrel\to\over {D \in Z^p(\ol{S})}}} 
\cl_{r,1}(\big|\ol{\rho}^{-1}(D)\big|) :
$$
$$
 \lim_{{\buildrel\to\over{D \in Z^p(\ol{S})}}} 
 \ch^r(\ol{\X}\bs \big|\ol{\rho}^{-1}(D)\big|,1;\Q) \to 
\lim_{{\buildrel\to\over{D \in Z^p(\ol{S})}}} \Gamma\big( 
H^{2r-1}(\ol{\X}\bs \big|\ol{\rho}^{-1}(D)\big|,\Q(r))\big).
$$

\proclaim{Proposition 3.4}
Assume the HC. Then,
$$
\coker (\cl^{p}_{r,1}) \cong \frac{\fL^p_{\ol{\X}/\ol{S}} 
\ch^r(\ol{\X};\Q)}{N_{\ol{S}}^{p} 
\ch^{r}(\ol{\X};\Q)}, \qquad \text{for $0 \leq p \leq r$.}
$$
Further, one does not need the HC assumption if 
$p\geq \text{\rm min}\{\dim\ol{\X} -3,r-1\}$.
\endproclaim

\demo{Proof} The proof proceeds exactly the same
way as in \S2 in obtaining the expression in (2.5),
with $D$ now a closed subscheme of $\ol{S}$ of (pure) codimension $p$,
and $\Y = \big|\ol{\rho}^{-1}(D)\big|$ being a closed 
subscheme of $\ol{\X}$ of codimension $p$, and where
$N^{1}_{\ol{S}}$ in (2.5) is replaced by $N^{p}_{\ol{S}}$. As in \S2, the
surjectivity of $\beta_{\Y}$ in diagram (2.3) is only guaranteed if 
the homological version of the HC (\cite{J3})
holds for $\Y$. For $p\geq \dim\ol{\X} -3$, that is clearly the
case since $\dim\Y \leq 3$. A similar story holds for
$p\geq r-1$, by the Lefschetz $(1,1)$ theorem.
\qed
\enddemo

\proclaim{Proposition 3.5}
Assume the HC. For any $p$ such that $0 \leq p \leq 
r,$ we have
$$
\fL^{p}_{\ol{\X}/\ol{S}} \ch^r(\ol{\X};\Q) = 
N_{\ol{S}}^p \ch^r(\ol{\X};\Q)\ +\ \ch^r(\ol{\X};\Q)_{AJ} \bigcap \K_p.
$$
Further, one does not need the HC assumption when
both conditions hold: $p\geq \text{\rm min}\{\dim\ol{\X} -3,r-1\}$ and 
$N_{\ol{S}}^{p}H^{2r-1}(\ol{\X},\QQ(r)) = 0$.
\endproclaim

\demo{Proof} Let 
$\xi \in \fL^{p}_{\ol{\X}/\ol{S}}\ch^r(\ol{\X};\Q))$.
Then by definition, there
is $\Y = \big|\ol{\rho}^{-1}(D)\big|\in Z^p(\ol{\X})$ 
such that if $j : \Y\hookrightarrow \ol{\X}$ is the inclusion,
then $\xi \in \IM(\a_{\Y})$ represents an element
in the numerator of the right hand side of the
following expression derived similarly to that in (2.4):
$$
\text{\rm coker}(\cl_{r,1}(\Y)) \simeq
\frac{\ker\biggl[ \ul{AJ}\big|_{\IM(\alpha_{\Y})}:
\IM(\alpha_{\Y}) \to J\biggl(
\frac{H^{2r-1}(\ol{\X},\QQ(r))}{j_{\ast}H_{\Y}^{2r-1}(\ol{\X},\QQ(r))}\biggr)
\biggr]}{\alpha_{\Y}\big(\ker(\beta_{\Y})\big)}.
$$
By the HC, the inclusion 
$j_{\ast}H^{2r-1}_{\Y}(\ol{\X},\QQ(r))
\hookrightarrow H^{2r-1}(\ol{\X},\QQ(r))$ has a cycle induced left inverse
$$
[w]_{\ast} : H^{2r-1}(\ol{\X},\QQ(r)) \twoheadrightarrow
j_{\ast}H^{2r-1}_{\Y}(\ol{\X},\QQ(r)) \subset H^{2r-1}(\ol{\X},\QQ(r)),
$$
where $|w|\subset \ol{\X}\times \Y$. 
Using functoriality of the Abel-Jacobi map and the fact that $\IM(\a_{\Y})
\subset \K_p$,  we have
$\xi - \alpha_{\Y}(w_{\ast}(\xi)) \in \ker (AJ) \bigcap \K_p$, and $w_{\ast}(\xi) \in
\ker(\beta_{\Y})$.
Thus
$$
\xi\in \ker(AJ)\bigcap \K_{p}\ +\ \alpha_{\Y}
( \ker(\beta_{\Y})) \subset \ker(AJ)\bigcap \K_p \ + \
N_{\ol{S}}^p\ch^r(\ol{\X};\Q).
$$
Hence this shows that
$$
\fL^{p}_{\ol{\X}/\ol{S}}\ch^r(\ol{\X};\Q) \subset \ker(AJ) \bigcap \K_p
+ N_{\ol{S}}^p\CH^r(\ol{\X};\QQ).
$$
Also by using Proposition 3.3, we have
$$
\fL^{p}_{\ol{\X}/\ol{S}}\ch^r(\ol{\X};\Q) \subset \ker(AJ)\bigcap \K_p\  + \
N_{\ol{S}}^p\CH^r(\ol{\X};\QQ) \subset 
\fL^{p}_{\ol{\X}/\ol{S}}\ch^r(\ol{\X};\Q)), \quad \text{for $1 \leq p \leq r$.} 
$$
\qed
\enddemo

\proclaim{Corollary 3.6}
  Assume the HC. Then for any $p$ such that $1 \leq p 
\leq r$, 
$$
\coker(\cl^p_{r,1}) \cong 
\frac{\fL^{p}_{\ol{\X}/\ol{S}}\ch^r(\ol{\X};\Q)}{N_{\ol{S}}^p
\ch^r(\ol{\X};\Q)}\cong 
\frac{\ch^r(\ol{\X};\Q)_{AJ}\bigcap \K_p}{\ch^r(\ol{\X};\Q)_{AJ} \bigcap \ker\K_p\bigcap
N_{\ol{S}}^p\ch^r(\ol{\X};\Q)}
$$
Further, the HC assumption is not needed if
both conditions hold:
$p\geq \text{\rm min}\{\dim\ol{\X} -3,r-1\}$ and 
$N_{\ol{S}}^{p}H^{2r-1}(\ol{\X},\QQ(r)) = 0$.
\endproclaim

For the remainder of this section, and indeed for the rest of this paper,
we are going to restrict to the special situation in Example 3.2 
with $X := \ol{\X} = \ol{S}$ and where
$\ol{\rho}$ is the identity.
\bigskip
\noindent
\remark{Remarks} 3.7. (i) Note that Corollary 3.6 implies
that
$$
\coker(\cl^r_{r,1}) \cong \ch^r(X;\Q)_{AJ}.
$$
Since the right hand side need not be
zero,  we recover the counterexample in \cite{J3}.
\bigskip
\noindent
(ii) Assuming the HC, 
together with $N^{r-1}\ch^r(X;\Q) = \ch^r_{\alg}(X;\Q)$
(\cite{J1} (Lemma 5.7)), and $J(N^{r-1}H^{2r-1}(X,\Q(r)))= 
J^r_a(X) := AJ(\ch^r_{\alg}(X;\Q))$ (see \cite{L}),
the reader can check the existence of the commutative diagram:
$$
\matrix
&&0 && 0 &&0&& \\
&\\
&& \downarrow && \downarrow&& \downarrow&& \\
&\\
0 &\to&  \ch^r_{\alg}(X;\Q)_{AJ} &\to& \ch^r(X;\Q)_{AJ} &\to& 
\frac{\fL^{r-1} \ch^r(X;\Q)}{N^{r-1}\ch^r(X;\Q)} &\to& 0\\
&\\
&& \downarrow && \downarrow && \downarrow
&& \\
&\\
 0 &\to& \ch^r_{\alg}(X;\Q) &\to& \ch^r_{\hom}(X;\Q) &\to &  
\Griff^r(X)_{\Q}&\to& 0 \\
&\\
 && \downarrow\rlap{\rm $AJ^{r-1}_a$} && \downarrow\rlap{\rm $AJ$} 
&\searrow\rlap{\rm$\ul{AJ}^{r-1}$} & \downarrow&& \\
&\\
 0 &\to& J^{r}_a(X) &\to & J^r(X) &\to& 
J\left(\frac{H^{2r-1}(X,\Q(r))}{N^{r-1}H^{2r-1}(X,\Q(r))} \right) 
&\to& 0\\
&\\
 && \downarrow && && && \\
 &\\
 && 0 &&&&&&
\endmatrix.\tag{3.8}
$$
where $AJ^{r-1}_a = AJ \big|_{\ch^r_{\alg}(X;\Q)}$, 
$\ch^r_{\alg}(X;\Q)_{AJ}=\ker(AJ^{r-1}_a)$. 
\endremark

%############ New section ###############
\head \S4. Level of Chow groups \endhead

We recall the following notion introduced in \cite{L}, that measures 
the complexity of Chow groups.

\proclaim{Definition 4.1} {\rm  We define
$$
\align
\level(H^{*}(X)) &= \text{\rm max}\big\{|p-q|\ \big|\ H^{p,q}(X)\ne 0,\ 
p,\ q \geq 0\big\},\\
\level(\ch^{r}(X;\Q)) &= \text{\rm min}\left\{\matrix &|&\ch^{r}(X\bs 
Y;\QQ) = 0\\
\ell\geq 0&|&\text{\rm over\ all\ closed}\ Y\subset X\\
&|&\text{\rm with}\ \codim_{X}Y = r-\ell
\endmatrix\right\},\\
\level(\ch^{*}(X;\Q)) &= \text{\rm 
max}\big\{\level(\ch^{r}(X;\Q))\ \big|\ r\geq 0\big\}.
\endalign
$$}
\endproclaim
\noindent
The notion of a conjectured descending Bloch-Beilinson (B-B) filtration
$\{F^{\nu}\ch^{r}(X;\Q)\}_{\nu=0}^{r}$  is widely introduced in the
literature. We refer the reader to \cite{J1} for the definition,
as well as the fact that if it exists, then it is unique. Among
the many properties of this filtration, one has
$F^{0}\ch^{r}(X;\Q) = \ch^{r}(X;\Q),\ F^{1}\ch^{r}(X;\Q)  = \ch^{r}_{\hom}
(X;\Q), \ F^{2}\ch^{r}(X;\Q) \subset \ker(AJ)$, and that the
induced action of a correspondence on $\text{\rm Gr}_{F}^{\nu}$ depends
only on the cohomology class of that correspondence. Various candidate 
B-B filtrations have been introduced in the literature, such as
by Griffiths/Green, J. P. Murre, M. Saito/M. Asakura, S. Saito, W. 
Raskind, J. D. Lewis, U. Jannsen, et al. {\it Under the aforementioned uniqueness 
property, we will make the blanket assumption that all 
such candidates define the same filtration.}

\proclaim{Theorem 4.2}
Assume GHC and existence of the B-B
filtration. Then for $\ell \geq 1$,
$$
\level(\ch^*(X;\Q)) \leq \ell \quad \Rightarrow \quad 
\ch^r_{\hom}(X;\Q) \subset N^{r-\ell} \ch^r(X;\Q) 
$$
and hence $\ch^r_{\hom}(X;\Q) = \fL^p\ch^r(X;\Q) = N^p\ch^r(X;\Q)$ 
for any $p$ satisfying $1 \leq p \leq r-\ell$.
 \endproclaim

\demo{Proof} 
First note that $\level(H^*(X)) = \level(\ch^*(X;\Q)) \leq \ell$ 
(\cite{L}, Corollary 15.64), hence 
$N^{p}H^{2r}(X,\Q(r))=H^{2r}(X,\Q(r))$ for any $p \leq (2r-\ell)/2$. 
By definition of level and by the GHC, there is a closed $Y \subset X$ 
possibly reducible of codimension $r-\ell$ such that both of 
$\ch^r_Y(X;\Q) \twoheadrightarrow \ch^r(X;\Q)$ and $H^{2r}_Y(X,\Q(r)) 
\twoheadrightarrow H^{2r}(X,\Q(r))=N^{r-\ell}H^{2r}(X,\Q(r))$ are 
surjective.  Let $\tilde{Y} \to Y$ be a desingularization of $Y $ and $\sigma: 
\tilde{Y} \to Y \hookrightarrow X$ be the composition of 
desingularization followed by an inclusion. Then we have the 
following diagram:
$$
\matrix
\ch^{\ell}(\tilde{Y};\Q) &@>\sigma_* >>& \ch^r(X;\Q) & \\
&\\
\downarrow\rlap{\rm $\cl^{\tY}_{\ell}$} & & \downarrow\rlap{\rm 
$\cl^X_r$}\\
&\\
H^{2 \ell}(\tilde{Y},\Q(\ell)) &@>[\sigma]_*>> & H^{2r}(X,\Q(r)) &= & 
N^{r-\ell} H^{2r}(X,\Q(r))
\endmatrix
$$
Note that $[\sigma]_*$ is surjective by Deligne's mixed Hodge theory 
(\cite{D} Proposition 8.2.8). Then by the HC, there is a 
cycle induced map 
$$
[\Gamma]_* : H^{2r}(X,\Q(r)) \to H^{2\ell}(\tilde{Y},\Q(\ell))
$$ 
such that $[\sigma]_* \circ [\Gamma]_* = \text{id}$ on 
$H^{2r}(X,\Q(r))$.  We now assume our given B-B filtration $F^{\bullet}\CH^r(-;\QQ)$,
let $\xi_1  \in F^{p} \ch^r(X;\Q)$, and set $\xi_0 := \Gamma_* (\xi_1) 
\in \ch^{\ell}(\tilde{Y};\Q)$. Then $\xi_0  \in F^{p} 
\ch^{\ell}(\tilde{Y};\Q)$ since the B-B filtration 
respects correspondences. 
Consider the following commutative diagram:
$$
\matrix
0 &\to & F^{p+1} \ch^{\ell} (\tY;\Q) & \to & F^{p}\ch^{\ell}(\tY;\Q)  
&@>{q^p_{\tY}}>> & \text{Gr}^{p}_F \ch^{\ell}(\tY;\Q) &\to&0 \\
&\\
&& \downarrow\llap{\rm$\sigma_*$\quad } \uparrow\rlap{\rm 
$\Gamma_*$}&&  \downarrow\llap{\rm$\sigma_*$\quad } \uparrow\rlap{\rm 
$\Gamma_*$} &&  \downarrow\llap{\rm$\sigma_*$\quad } 
\uparrow\rlap{\rm $\Gamma_*$}&& \\
&\\
0 &\to& F^{p+1} \ch^{r} (X;\Q) &\to & F^{p}\ch^{r}(X;\Q) 
&@>{q^{p}_{X}}>>  & \text{Gr}^{p}_F \ch^{r}(X;\Q) &\to & 0
\endmatrix
$$
\medbreak
\noindent
\underbar{Claim 1}: $q^{p}_X (\xi_1 - \sigma_*(\xi_0))=0$ :  
\medbreak
\noindent
Note that on graded level $\text{Gr}^{\dt}_F$, the action of 
correspondences depends only on the cohomology class of such 
correspondences, so
$$
\align
  q^{p}_X (\xi_1 - \sigma_*(\xi_0))&= q^{p}_X(\xi_1) - 
q^{p}_X(\sigma_*(\Gamma_*(\xi_1)) =  q^{p}_X(\xi_1) -\sigma_* \circ 
q^{p}_{\tY}(\Gamma_*(\xi_1))\\
&=  q^{p}_X(\xi_1) - \sigma_*(\Gamma_* \circ q^{p}_X (\xi_1))= 
q^{p}_X(\xi_1) -  ([\sigma]_* \circ [\Gamma]_* )(q^{p}_X(\xi_1))=0.
\endalign
$$
\medbreak
\noindent
\underbar{Claim 2}: For any $\xi \in \ch^r_{\hom}(X;\Q)= F^1 \ch^r(X;\Q)$, 
there is $\xi_Y \in F^1 \ch^{\ell}(\tY;\Q)$  and $\zeta_{\ell} \in 
F^{\ell+1}\ch^r(X;\Q) $ such that $\xi_1 = 
\sigma_*(\xi_Y)+\zeta_{\ell}$ :
\medbreak
\noindent
Let $\xi \in \ch^r_{\hom}(X;\Q)=F^1 \ch^r(X;\Q)$. Set $\xi_1 = 
\Gamma_* (\xi) \in F^1 \ch^{\ell}(\tY;\Q)$. Then by Claim 1, 
$$
\zeta_1 := \xi - \sigma_* (\xi_1) \in \ker (q^1_X)=F^2 \ch^r(X;\Q)
$$
Set $\xi_2 = \Gamma_*(\zeta_1) \in F^2\ch^{\ell}(\tY;\Q)$. By 
applying Claim 1 to $\zeta_1 \in F^2 \ch^r(X;\Q)$, we get 
$$
\zeta_2:= \zeta_1 - \sigma_*(\xi_2) \in \ker (q^2_{X}) = F^3 
\ch^r(X;\Q)
$$
By repeating this process $\ell$ times, we arrive at
$$
\xi=\sigma_*(\xi_1) + \sigma_* (\xi_2) + \cdots +\sigma_*(\xi_{\ell}) 
+ \zeta_{\ell},
$$
where $\xi_i \in F^i\ch^{\ell}(\tY;\Q)$ and $\zeta_{\ell} \in 
F^{\ell+1}\ch^{r}(X;\Q)$. Set $\xi_Y =\sum^{\ell}_{i=1} \xi_i \in F^1 
\ch^{\ell}(\tY;\Q)$.
\medbreak
\noindent
\underbar{Claim 3}: $F^{\ell+1} \ch^r(X;\Q)=0$ : 
\medbreak
\noindent
Note that
$$
\text{Gr}_{F}^{p}\ch^{r}(X;\Q) = \Delta_{X}(2d-2r+p,2r-p)_{\ast}
\ch^{r}(X;\Q).
$$
where $d= \dim X$ and $\Delta_X \in \ch^d(X\times X ;\Q)$ is the 
diagonal class.
Then,  
$$
[\Delta_{X}(2d-2r+p,2r-p)] \in H^{2d-2r+p}(X,\Q)
\otimes H^{2r-p}(X,\Q),
$$
and by the GHC together with $\level(H^{\ast}(X))\leq \ell$
we can assume that $[\Delta_{X}(2d-2r+p,2r-p)]$ is
supported on some $W_{1,p}\times W_{p,2}\subset X\times X$ with 
$\codim(W_{1,p}, X) \geq (2d-2r+p-\ell)/2$ and 
$\codim(W_{p,2},X) \geq (2r-p-\ell)/2$, where one can 
easily check that
$\Delta_{X}(2d-2r+p,2r-p)_{\ast} = 0$ on
$\ch^{r}(X;\Q)$ for $p \geq \ell+1$. Hence 
$F^{\ell+1}\ch^r(X;\Q) =  0$.
\medbreak
\noindent
From Claims 1, 2 and 3, we get, for any $\xi \in \ch^r_{\hom}(X;\Q)$, 
$$
\xi = \sigma_* \left( \sum^{\ell}_{i=1} \xi_i\right) \in \sigma_* 
\left( F^1\ch^{\ell}(\tY;\Q) \right)
$$
since $F^{\dt}\ch^{\ell}(\tY;\Q)$ is a descending filtration. 
\medbreak
\noindent
\underbar{Claim 4}: $\sigma_* (F^1 \ch^{\ell}(\tY;\Q)) \subset 
N^{r-\ell}\ch^r(X;\Q)$ :
\medbreak
\noindent
This follows easily from the fact that $F^1 \ch^{\ell}(\tY;\Q) = \ch^{\ell}_{\hom}(\tY;\Q)$.
\medskip
\noindent
By the above claims, we have now completed the proof of the lemma. \qed
 \enddemo

Recall that $N^{r-1}\CH^r(X;\QQ) = \CH^r_{\alg}(X;\QQ)$. We deduce:
\proclaim{Corollary 4.3}
Assume the GHC and the existence of the B-B filtration. Then
$$
\level(\CH^{\ast}(X;\QQ)) \leq 1 \Rightarrow
\Griff^{\ast}(X)_{\QQ}= 0.
$$
\endproclaim

\remark{Remark 4.4}
The converse of Corollary 4.3 is false, as can be seen
by taking $X$ to be any surface of positive
geometric genus. 
\endremark

%################### New section 
%#############################################

 \head \S5. Summary Consequences \endhead

\proclaim{Proposition 5.1}
Assume the HC and a given $p$. Then the following statements are 
equivalent: 
\medbreak\noindent
{\rm (i)} $\cl^{p}_{r,1}$ is surjective,
\medbreak\noindent
{\rm (ii)} $\ch^r(X;\Q)_{AJ} \subset N^p \ch^r(X;\Q)$,
\medbreak\noindent
{\rm (iii)} $\cl^{\nu}_{r,1}$ is surjective for all $1 \leq \nu 
\leq p$.
\endproclaim

\demo{Proof} Use Corollary 3.6.\qed
\enddemo

\proclaim{Corollary 5.2}
Assume the GHC and existence of the  B-B filtration. Then for 
$\ell \geq 1$, if $\level(\ch^*(X;\QQ)) \leq \ell$, then 
$\cl^p_{r,1}$ is surjective for all $p$ satisfying $1 \leq p \leq 
r-\ell$. 
\endproclaim

\demo{Proof} 
In order to show the surjectivity of $\cl^p_{r,1}$ for $1 \leq p 
\leq r-\ell$, it is enough by Proposition 5.1 
to show that $\cl^{r-\ell}_{r,1}$ is 
surjective, equivalently 
$$
\coker(\cl^{r-\ell}_{r,1}) \cong {\fL^{r-\ell}\ch^r(X;\Q) \over 
N^{r-\ell}\ch^r(X;\Q)}=0
$$
First note that $\level(H^*(X)) \leq \ell$ implies $N^{r - \ell} 
H^{2r-1}(X,\Q(r))= H^{2r-1}(X,\Q(r))$ and in turn
$$
\align
\fL^{r - \ell}\ch^r(X;\Q) &=\ker  \left[\ch^r_{\hom}(X;\Q) \to 
J\biggl(\frac{H^{2r-1}(X,\Q)}{N^{r- \ell}H^{2r-1}(X,\Q(r))} 
\biggr)\right] \\
&= \ch^r_{\hom}(X;\Q) \subset N^{r-\ell}\ch^r(X;\Q),
\endalign
$$
where the last inclusion follows from Theorem 4.2.  Hence we are done. 
\qed
\enddemo

In case when $\ell=1$, Corollary 5.2 gives the surjectivity of  
$\cl^p_{r,1}$ up to $p =r-1$. However we can extend this to the 
surjectivity of $\cl^r_{r,1}$ by the result of S. Saito 
$F^2\ch^r_{\alg}(X;\Q)= \ch^r(X;\Q)_{AJ} \cap 
\ch^r_{\alg}(X;\Q)$ (\cite{S} Corollary 3.7). In fact, we have 
following:

\proclaim{Corollary 5.3}
Assume the GHC and existence of the  B-B filtration. Then the 
following statements are equivalent:
  \medbreak\noindent
{\rm (i)} $\level (\ch^*(X;\Q)) \leq 1$,
\medbreak\noindent
{\rm (ii)} $\cl^p_{r,1}$ is surjective for any $r,\ p$,
\medbreak\noindent
{\rm (iii)} $\ch^*_{\alg}(X;\Q) \cong J^*_a (X)$,
\medbreak\noindent
{\rm (iv)} $\ch^*_{\alg}(X;\Q)$ is representable.
\endproclaim
 
\demo{Proof}
${\text {\rm (ii)}} \Rightarrow {\text {\rm (iii)}}$ : 
Suppose $\cl^p_{r,1}$ is surjective for 
all $p,r$. Then, in particular, $\cl^r_{r,1}$ is surjective for 
any $r$ and hence $\ch^r(X;\Q)_{AJ}=0$ for any $r$ (Example 3.2). 
This implies that  $AJ^{r-1}_a: \ch^r_{\alg}(X;\Q) @>\simeq>> J^r_a(X)$ is an 
isomorphism for any $r$. Thus $\ch^*_{\alg}(X;\Q) \cong J^*_a(X;\Q)$.
\medbreak
\noindent
${\text {\rm (iii)}} \Rightarrow {\text {\rm (iv)}}$ : Obvious.
\medbreak
\noindent
${\text {\rm (iv)}} \Rightarrow {\text {\rm (i)}}$ : \cite{L} (Corollary 15.42).
\medbreak
\noindent
${\text {\rm (i)}} \Rightarrow {\text {\rm (ii)}}$ : 
Suppose $\level (\ch^*(X;\QQ))\leq 1$. 
We refer to diagram (3.8). Then by Corollary 4.3,
$$
\frac{\fL^{r-1}\ch^r(X;\Q)}{N^{r-1}\ch^r(X;\Q)} = 0,\ \text{\rm hence}\
\ch^{r}_{\alg}(X;\Q)_{AJ} = \ch^r(X;\Q)_{AJ}.
$$ 
Now statement (ii) follows from  Claim 3 in 
the proof of Theorem 4.2, since in this
case $\level (\ch^*(X;\Q)) \leq 1
\Rightarrow F^{2}\CH^{r}(X;\QQ) =0$, together with
the aforementioned result due to S. Saito, viz.,
$\ch^r_{\alg}(X;\Q)_{AJ} = F^2 \ch^r_{\alg}(X;\Q)$. \qed
\enddemo

%################   New section  
%##########################################

\head \S6. Some key examples \endhead

Our first example, which illustrates the nontriviality of 
coker$(\cl_{r,1}^{\nu})$ for $2\leq\nu\leq r$,
is based on a rewording of a theorem in Nori's paper \cite{N}.
In order to reword Nori's theorem, we have to make the following 
assumptions:
\medbreak\noindent
(1) GHC
\medbreak\noindent
(2) The conjecture in Jannsen's paper (\cite{J1}) about Nori's
filtration. Specifically  
$$
A_{\ell}\ch^{r}(X;\Q) = 
N^{\nu}\ch^{r}(X;\Q), \quad\nu = r-\ell-1,
$$
where $A_{\ell}\ch^{r}(X;\Q)$ is defined in \cite{N}.
Nori's theorem translated in our language is now the following:
\proclaim{Theorem 6.0} {\rm (\cite{N})}
Let $W$ be a smooth complex projective variety
and $X \subset W$ a sufficiently general complete intersection of
sufficiently high multidegree. Let $\xi_{W}\in \ch^{r}(W;\Q)$
and put $\xi := \xi_{W}\big|_{X}$. Assume $2r-\nu -1 < d := \dim X$. 
If 
$\xi \in N^{\nu}\ch^{r}(X;\Q)$, then:
\medbreak\noindent
{\rm (1)} $\big[\xi_{W}\big] = 0\in H^{2r}(W,\Q)$,
\medbreak\noindent
{\rm (2)} $\underline{AJ}^{\nu}(\xi_{W}) = 0
\in 
\text{\rm Ext}^{1}_{\text{MHS}}\biggl(\Q(0),\frac{H^{2r-1}(W,\Q(r))}{N^{\nu}
H^{2r-1}(W,\Q(r))}\biggr)$.
\endproclaim

\noindent{\underbar{Example 6.1}.} The following example is essentially due 
to Nori (\cite{N}). Let
$W = Q_{2r}\subset \PP^{2r+1}$ be a smooth quadric. Let $\xi_W\in 
\ch^r(W;\Q)$
be given such that $\text{Prim}^{2r}(W;\Q) = \Q\big[\xi_W\big]$. Let 
$X\subset W$ be
a general complete intersection of sufficiently high multidegree of 
dimension $d\leq 2r-1$.
Note that $H^{2r-1}(W,\Q) = 0$, and we also want 
to arrange for $H^{2r-1}(X,\Q) = 0$. 
%(Further, note that since
%$\big[\xi_W\big] \ne 0$,  $AJ(\xi_W)$ is undefined.) 
To ensure that 
$H^{2r-1}(X,\Q) = 0$, we must have $d \ne 2r-1$. In particular,
in the relation $2r-\nu-1 < d< 2r-1$, we must necessarily have 
$2\leq \nu\leq r$.  Put $\xi = \xi_W\big|_{X} \in \CH^r_{\hom}(X;\QQ)$. Then by Nori's theorem
with the range of $\nu$, 
$$
\xi \ne 0\in \frac{\fL^{\nu}\ch^r(X;\Q)}{N^{\nu}\ch^r(X;\Q)}.
$$
Thus coker$(\cl_{r,1}^{\nu}) \ne 0$ for $2\leq\nu\leq r$. Complementary results appear
in our next example.

\medskip
\noindent
{\underbar{Example 6.2}}. Suppose $X \subset \PP^{N}$ be a smooth complete 
intersection of dimension $d$. Then by using the Lefschetz theorem, 
for $\nu \geq 1$, $\text{Gr}^{\nu}_F \ch^r(X;\Q) = 0$ if $2r-\nu 
\neq d$ (S. Saito). Let $F^{\dt}\ch^r(X;\Q)$ be a the (unique)
conjectured B-B  filtration. It is known that 
$$
F^2\ch^r(X;\Q) \subset \ch^r(X;\Q)_{AJ} \subset F^1 \ch^r(X;\Q),
$$
and let us assume the following conjectural statement of Jannsen \cite{J1}:
$$
F^{\nu+1} \ch^r(X;\Q) \subset N^{\nu} \ch^r(X;\Q).
$$
For a given $p$, if we impose the condition that $\text{Gr}^{\nu}_F 
\ch^r(X;\Q) =0$ for all $1 \leq \nu \leq p$, (i.e. $F^1 \ch^r(X;\Q) = 
\cdots = F^{p+1}\ch^r(X;\Q)$), then this implies
$$
\ch^r(X;\Q)_{AJ} \subset F^1 \ch^r(X;\Q) = F^{p+1}\ch^r(X;\Q) \subset 
N^p\ch^r(X;\Q),
$$
hence by Proposition 5.1, $\cl^{\nu}_{r,1}$ is surjective for all $1 
\leq \nu \leq p$. Note that the requirement $2r-\nu \neq d$ for
all $1\leq \nu\leq p$ is equivalent to  $2r-d \not\in
\{1,\ldots,p\}$. Hence we have either $2r-d \leq 0$ or $2r-d \geq 
p+1$. 
\medbreak\noindent
{\rm (1)} If $2r-d \leq 0$, then for any $\nu \geq 1$, we have
$$
2r-\nu \leq d-\nu < d.
$$
Hence for any $\nu \geq 1$, $\cl^{\nu}_{r,1}$ is surjective.
\medbreak\noindent
{\rm (2)} If $2r-d \geq p+1$, then:
$$
\text{Gr}^{\nu}_F \ch^r(X;\Q) =0, \qquad \text{for $1 \leq \nu \leq 
2r-d-1$},
$$
and hence $\cl^{\nu}_{r,1}$ is surjective for $1 \leq \nu \leq 
2r-d-1$.

\remark{Remark 6.3}
We want to make it clear that in the previous
example, we have either $d> 2r-1$ or $d \leq 2r-\nu-1$. This is in total
contrast to  Example 6.1,
where $2r-\nu-1 < d < 2r-1$, and where in this case
$\cl_{r,1}^{\nu}$ is not surjective for $2\leq \nu\leq r$.
\endremark

\Refs

\item{[As]} M. Asakura, Motives and algebraic de Rham cohomology,
in: The Arithmetic and Geometry of Algebraic Cycles, Proceedings of 
the CRM Summer School, June 7-19, 1998, Banff, Alberta, Canada 
(Editors: B. Gordon, J. Lewis, S. M\"uller-Stach, S. Saito and N. 
Yui), CRM Proceedings and Lecture Notes {\bf 24}, AMS 2000, 133--155

\item{[Be]} A. Beilinson, Notes on absolute Hodge cohomology,
Contemporary Mathematics {\bf 55}, Part I, 1986, 35--68.

\item{[Blo]} S. Bloch, Algebraic cycles and higher $K$-theory,
Adv. Math {\bf 61}, (1986), 267--304.

\item{[C]} J. Carlson, Extensions of mixed Hodge structures,  
JournŽes de GŽometrie AlgŽbrique d'Angers, Juillet 1979/Algebraic 
Geometry, Angers, 1979,  pp. 107--127

\item{[CL]} Xi Chen, J. D. Lewis: The real regulator for a product of
$K3$ surfaces, Mirror Symmetry V, Proc. of the BIRS
Workshop on Calabi-Yau Varieties and Mirror Symmetry, 
Edited by N. Yui, S.-T. Yau, J. D.
Lewis, AMS/IP Studies in Advanced Mathematics, Vol.
{\bf 28}, (2006), 271--283.

%\item{[DL]} R. de Jeu, J. D. Lewis, Belinson's Hodge conjecture
%or smooth varieties. In preparation.

\item{[D]} P. Deligne, Theorie de Hodge III, Inst. Hautes \'Etudes Sci. 
Publ. Math. No. {\bf 44} (1974), 5--77. 

\item{[G]} A. Grothendieck, Hodge's general conjecture is false for 
trivial reasons,  Topology {\bf 8},  1969, 299--303.

\item{[J1]} U. Jannsen, Equivalence relations on algebraic cycles, in: The 
Arithmetic and Geometry of Algebraic Cycles, Proceedings of the NATO 
Advanced Study Institute on The Arithmetic and Geometry of Algebraic 
Cycles, Banff, Alberta, Canada, 1998, (Editors: B. Gordon, 
J. Lewis, S. M\"uller-Stach, S. Saito and N. Yui), Kluwer Academic 
Publishers, Dordrecht, The Netherlands {\bf 548}, (2000),
225--260.

\item{[J2]} U. Jannsen, Deligne cohomology, Hodge-${\Cal 
D}$-conjecture, and motives: in Beilinson's Conjectures on Special 
Values of $L$-Functions, (Editors: Rapoport, Schappacher, Sneider), 
Perspectives in Math. {\bf 4}, Academic Press (1988), 305--372.

\item{[J3]} U. Jannsen, Mixed Motives and Algebraic $K$-Theory, Lecture
Notes in Mathematics {\bf 1400}, (1988).

\item{[KL]} M. Kerr, J. D. Lewis, The Abel-Jacobi map for higher Chow
groups II, Inventiones Math. {\bf 170} (2), (2007), 355--420.

\item{[L]} J. D. Lewis,  A survey of the Hodge conjecture, Second 
edition. Appendix B by B. Brent Gordon. CRM Monograph Series, {\bf 10}. 
American Mathematical Society, Providence, RI, 1999.

%\item{[L2]} J. D. Lewis,  Cylinder homomorphisms and Chow groups,  
%Math. Nachr.  {\bf 160}  (1993), 205--221. 

%\item{[LS]} J. D. Lewis,  A. S. Sert\"{o}z, Motives of Some Fano 
%Varieties, preprint

\item{[N]} M. Nori,  Algebraic cycles and Hodge-theoretic 
connectivity,  Invent. Math.  111  (1993),  no. 2, 349--373.

\item{[S]} S. Saito, Motives and filtrations on Chow groups,  Invent. 
Math. {\bf 125}  (1996),  no. 1, 149--196.

\endRefs

\enddocument
\bye